\documentclass[11pt]{amsart}
\usepackage{mathrsfs}
\usepackage{amssymb}

\usepackage{titletoc}
\usepackage{amscd}
\usepackage{amsmath, amssymb}
\usepackage{amsfonts}
\usepackage[colorlinks,linkcolor=blue,citecolor=blue, pdfstartview=FitH]{hyperref}

\numberwithin{equation}{section}

\theoremstyle{plain}
\newtheorem{thm}{Theorem}[section]
\newtheorem{theorem}[thm]{Theorem}
\newtheorem{lemma}[thm]{Lemma}
\newtheorem{corollary}[thm]{Corollary}
\newtheorem{proposition}[thm]{Proposition}

\theoremstyle{definition}

\newtheorem{remark}[thm]{Remark}

\newtheorem{definition}[thm]{Definition}

\newtheorem{example}[thm]{Example}

\newtheorem{defn-thm}[thm]{Definition-Theorem}

\newcommand{\sA}{{\mathcal A}}

\newcommand{\sE}{{\mathcal E}}

\newcommand{\sL}{{\mathcal L}}
\newcommand{\sM}{{\mathcal M}}

\newcommand{\sR}{{\mathcal R}}

\newcommand{\sT}{{\mathcal T}}

\newcommand{\sX}{{\mathcal X}}

\newcommand{\B}{{\mathbb B}}
\newcommand{\C}{{\mathbb C}}

\renewcommand{\H}{{\mathbb H}}

\renewcommand{\P}{{\mathbb P}}

\newcommand{\R}{{\mathbb R}}

\newcommand{\qtq}[1]{\quad\mbox{#1}\quad}
\newcommand{\bp}{\bar{\partial}}
\newcommand{\Om}{\Omega}

\newcommand{\ts}{\otimes}

\newcommand{\btheorem}{\begin{theorem}}
\newcommand{\etheorem}{\end{theorem}}
\newcommand{\bproposition}{\begin{proposition}}
\newcommand{\eproposition}{\end{proposition}}
\newcommand{\bdefinition}{\begin{definition}}
\newcommand{\edefinition}{\end{definition}}
\newcommand{\bcorollary}{\begin{corollary}}
\newcommand{\ecorollary}{\end{corollary}}
\newcommand{\bproof}{\begin{proof}}
\newcommand{\eproof}{\end{proof}}
\newcommand{\bremark}{\begin{remark}}
\newcommand{\eremark}{\end{remark}}
\newcommand{\eexample}{\end{example}}
\newcommand{\bexample}{\begin{example}}
\newcommand{\la}{\langle}
\newcommand{\elemma}{\end{lemma}}
\newcommand{\blemma}{\begin{lemma}}
\newcommand{\ra}{\rangle}
\newcommand{\sq}{\sqrt{-1}}

\newcommand{\p}{\partial}

\renewcommand{\bar}{\overline}

\renewcommand{\phi}{\varphi}

\newcommand{\ee}{\end{eqnarray*}}
\newcommand{\be}{\begin{eqnarray*}}

\newcommand{\beq}{\begin{equation}}
\newcommand{\eeq}{\end{equation}}

\newcommand{\bd}{\begin{enumerate}}
\newcommand{\ed}{\end{enumerate}}

\renewcommand{\bf}{\textbf}

\renewcommand{\bf}{\textbf}

\renewcommand{\>}{\rightarrow}

\usepackage{fancyhdr}
\renewcommand{\leftmark}{{\scriptsize Curvatures of moduli space of curves and applications}}

\begin{document}
\title{Curvatures of moduli space of curves and applications}
\makeatletter
\let\uppercasenonmath\@gobble
\let\MakeUppercase\relax
\let\scshape\relax
\makeatother
\author{Kefeng Liu, Xiaofeng Sun, Xiaokui Yang and Shing-Tung Yau}

\address{Kefeng Liu, Department of Mathematics, UCLA}
\email{liu@math.ucla.edu}
\address{Xiaofeng Sun, Department of Mathematics, Lehigh University}
\email{xis205@lehigh.edu}

\address{Xiaokui Yang, Morningside Center of Mathematics, Academy of Mathematics and Systems Science, Chinese Academy of Sciences, Beijing, China;
 Hua Loo-Keng Key Laboratory of Mathematics, Academy of Mathematics and Systems Science, Chinese Academy of Sciences, Beijing, China} \email{xkyang@amss.ac.cn}

\address{Shing-Tung Yau, Department of Mathematics, Harvard University.}
\email{yau@math.harvard.edu}
 \maketitle

\begin{abstract}  In this paper, we investigate the geometry of the moduli
space of curves by using the curvature properties of direct image
sheaves of vector bundles. We show that the moduli space $(\sM_g,
\omega_{\mathrm{WP}})$ of curves with genus $g>1$ has dual-Nakano
negative and semi-Nakano-negative curvature, and in particular, it
has non-positive Riemannian curvature operator and also non-positive
complex sectional curvature. We also prove that any submanifold  in
$\sM_g$ which is totally geodesic in $\sA_g$ with finite volume must
be a ball quotient.

\end{abstract}

\section{Introduction}

In this paper, we study the curvature properties of the
Weil-Petersson metric as well as its background Riemannian metric on
the moduli space  of curves.

On Riemannian manifolds, there are many curvature terminologies,
e.g.  curvature operator, sectional curvature, isotropic curvature
and etc.. As it is well-known, the curvature relations are well
understood on Riemannian manifolds (e.g. \cite[p.100]{Brendle11}).
On the other hand, we also have some classical curvature concepts on
K\"ahler manifolds, such as the holomorphic bisectional curvature,
curvature in the sense of Siu and curvature in the sense of Nakano.
At first, we obtain a  list of curvature relations between a variety
of curvature properties of the K\"ahler metric and its background
Riemannian metric:

\btheorem\label{1} On a  K\"ahler manifold $(X,\omega)$, the
curvatures

\bd

\item semi dual-Nakano-negative;

\item non-positive Riemannian curvature operator;

\item strongly non-positive in the sense of siu;

\item non-positive complex sectional curvature;

\item non-positive Riemannian sectional curvature;

\item non-positive holomorphic bisectional curvature;

\item non-positive isotropic curvature
\ed have the following relations
$$ (1)\Longrightarrow(2)\Longrightarrow (3)\Longleftrightarrow(4)\Longrightarrow (5)\Longrightarrow (6);$$
$$(1)\Longrightarrow (3)\Longleftrightarrow(4)\Longrightarrow (7).$$
\etheorem

Let $f:\sT_g\>\sM_g$ be  the universal curve with genus $g\ge 2$.
Since it is a canonically polarized family,  and
$$T^{*}\sM_g\cong f_*(K^{\ts 2}_{\sT_g/\sM_g}), $$
one can compute the curvature of the induced metric on $T^*\sM_g$ by
using the curvature formula of direct image sheaves (e.g.
\cite{Berndtsson11}, \cite{Schumacher12} and \cite{Liu-Yang}). This
induced metric is actually  conjugate dual to the Weil-Petersson
metric $\omega_{\mathrm{WP}}$ on $\sM_g$. By adapting the methods in
\cite{Liu-Sun-Yau08}, we show that $(\sM_g,\omega_{\mathrm{WP}})$
has the similar curvature properties as the space form--the unit
disk $(\B^{3g-3},\omega_{\mathrm{B}})$ with the invariant Bergman
metric $\omega_{\mathrm{B}}$, i.e. $(\sM_g,\omega_{\mathrm{WP}})$
possesses the strongest curvature properties of complex manifolds.

\btheorem\label{2} The curvature of Weil-Petersson metric
$\omega_{\mathrm{WP}}$ on the moduli space $\sM_g$ of Riemann
surfaces of genus $g\geq 2$ is dual-Nakano-negative and
semi-Nakano-negative.  \etheorem

\noindent As applications of Theorem \ref{1} and Theorem \ref{2}, we
obtain a variety of curvature properties of the moduli space of
curves:

\btheorem\label{3} The moduli space $(\sM_g, \omega_{\mathrm{WP}})$
 has the following curvature properties:

\bd

\item dual-Nakano-negative and  semi Nakano-negative curvature;

\item non-positive Riemannian  curvature operator;

\item non-positive complex sectional curvature;

\item strongly-negative curvature in the sense of Siu;

\item negative Riemannian sectional curvature;

\item negative holomorphic bisectional curvature;

\item non-positive isotropic curvature.

\ed \etheorem

\noindent Note that, part  $(2)$ is firstly obtained in \cite{Wu}
recently.

We  describe another application of the curvature properties of
$(\sM_g, \omega_{\mathrm{WP}})$. Denote by $\sA_{g}$ the moduli
space of principally polarized abelian varieties of dimension $g$
and denote by \beq j: \sM_g\>{\sA_g}\eeq the Torelli map associating
to a curve its Jacobian with its natural principal polarization.  We
denote by $Jac_g$  the image $j(\sM_g)$ and let $\bar{Jac_g}$, the
so-called
 Torelli locus, be the schematic closure of $Jac_g$ in $\sA_g$. Frans
Oort asked in \cite[Section~7]{Oort95} whether there exists any
locally symmetric subvariety of $\sA_g$  which is contained in
${\bar {Jac_g}}$ and intersects $Jac_g$, and he conjectured that
nontrivial such subvarieties do not exist. This conjecture is
extensively studied in the last decade by using algebraic geometry
methods. In this paper, we use a differential geometric approach and
obtain the following

\btheorem\label{4} Let $j:(\sM_g,
\omega_{\mathrm{WP}})\>(\sA_g,\omega_{H})$ be the Torelli map  where
$\omega_H$ is the Hodge metric. Let $V$ be a submanifold in $\sM_g$
with
 $j(V)$  totally geodesic in $(\sA_g,
\omega_{H})$. If $j(V)$ has finite volume, then $V$ must be a ball
quotient. In particular, any compact submanifold $V$ in $\sM_g$ with
 $j(V)$  totally geodesic in $(\sA_g,\omega_{H})$ must be a ball
 quotient.

 \etheorem

\bremark  By using algebraic  methods, Hain
(\cite[Theorem~1]{Hain99}) and de Jong-Zhang
(\cite[Theorem~1.1]{deJ-Zhang07}) proved similar results under
certain conditions. See also \cite{MVZ05}, \cite{Lu-Zuo} and
\cite{CFG}. For more progress on Oort's conjecture, we refer the
reader to   survey papers  \cite{Oort95, Moonen-Oort}.

\eremark

\noindent As a special case, we show that there is no higher rank
locally symmetric space in $\sM_g$:

\bcorollary\label{5} Let $\Om$ be an irreducible bounded symmetric
domain and $\Gamma\subset Aut(X)$ be a torsion-free cocompact
lattice, $X := \Om/\Gamma$. Let $h$ be the canonical metric on $X$.
If there exists a  nonconstant holomorphic mapping $f:(X,h)\>(\sM_g,
\omega_{\mathrm{WP}})$, then $\Om$ must be of rank $1$, i.e. $X$
must be a ball quotient.
 \ecorollary

\bf{Acknowledgement}. The authors would  like to thank  the
anonymous referee for clarifying many issues of the paper. The third
named author would also like to thank Valentino Tosatti for many
helpful suggestions.

\section{Curvature relations on  vector bundles}

Let $E$ be a holomorphic vector bundle over a  K\"ahler manifold $X$
and $h$ a Hermitian metric on $E$. There exists a unique connection
$\nabla$ which is compatible with the
 metric $h$ and the complex structure on $E$. It is called the Chern connection of $(E,h)$. Let $\{z^i\}_{i=1}^n$ be  the local holomorphic coordinates
  on $X$ and  $\{e_\alpha\}_{\alpha=1}^r$ be a local frame
 of $E$. The curvature tensor $R^\nabla\in \Gamma(X,\Lambda^2T^*X\ts E^*\ts E)$ has the form
 \beq R^\nabla=\frac{\sq}{2\pi} R_{i\bar j\alpha}^\gamma dz^i\wedge d\bar z^j\ts e^\alpha\ts e_\gamma\eeq
where $R_{i\bar j\alpha}^\gamma=h^{\gamma\bar\beta}R_{i\bar
j\alpha\bar \beta}$ and \beq R_{i\bar j\alpha\bar\beta}= -\frac{\p^2
h_{\alpha\bar \beta}}{\p z^i\p\bar z^j}+h^{\gamma\bar
\delta}\frac{\p h_{\alpha \bar \delta}}{\p z^i}\frac{\p
h_{\gamma\bar\beta}}{\p \bar z^j}.\label{dcurvature}\eeq Here and
henceforth we
 adopt the Einstein convention for summation.

\bdefinition
 A Hermitian vector bundle
$(E,h)$ is said to be \emph{Griffiths-positive}, if for any nonzero
vectors $u=u^i\frac{\p}{\p z^i}$ and $v=v^\alpha e_\alpha$,  \beq
\sum_{i,j,\alpha,\beta}R_{i\bar j\alpha\bar \beta}u^i\bar u^j
v^\alpha\bar v^\beta>0.\eeq $(E,h)$ is said to be
\emph{Nakano-positive}, if for any nonzero vector
$u=u^{i\alpha}\frac{\p}{\p z^i}\ts e_\alpha$, \beq
\sum_{i,j,\alpha,\beta}R_{i\bar j\alpha\bar \beta} u^{i\alpha}\bar
u^{j\beta}>0. \label{Nakanop}\eeq $(E,h)$ is said to be
\emph{dual-Nakano-positive}, if for any nonzero vector
$u=u^{i\alpha}\frac{\p}{\p z^i}\ts e_\alpha$, \beq
\sum_{i,j,\alpha,\beta}R_{i\bar j\alpha\bar \beta} u^{i\beta}\bar
u^{j\alpha}>0. \label{dual}\eeq \noindent The notions of
semi-positivity, negativity and semi-negativity can be defined
similarly. We say $E$ is Nakano-positive (resp. Griffiths-positive,
dual-Nakano-positive, $\cdots$), if it admits a
Nakano-positive(resp. Griffiths-positive, dual-Nakano-positive,
$\cdots$) metric. \edefinition

\bremark\label{dual123} It is easy to see that $(E,h)$ is
dual-Nakano-positive if and only if $(E^*,h^*)$ is Nakano-negative.
\eremark

As  models of complex manifolds, one has the following well-known
curvature properties:

\blemma\label{optimal} Let $n>1$.
 \bd\item $(T\P^n,\omega_{FS})$ is dual-Nakano-positive and
 semi-Nakano-positive.

\item Let $X$ be a  hyperbolic space form with dimension $n$. If $\omega_B$ is the canonical metric on
$X$, then $(TX,\omega_B)$ is dual-Nakano-negative and
semi-Nakano-negative.

 \ed
 \elemma

We shall use the following curvature monotonicity formulas
frequently, in particular the explicit curvature formulas(e.g.
(\ref{cu})). Hence we include a detailed proof.
\blemma\label{monotone} Let $(E,h)$ be a Hermitian holomorphic
vector bundle over a complex manifold $X$, $S$ be a holomorphic
subbudle of $E$ and $Q$ the corresponding quotient bundle,
$0\>S\>E\>Q\>0$.

\bd \item If $E$ is (semi-)Nakano-negative, then $S$ is also
(semi-)Nakano negative.

\item If $E$ is (semi-)dual-Nakano-positive, then $Q$ is also
(semi-)dual-Nakano-positive. \ed \bproof This lemma is
well-known(e.g.\cite{Demailly}).
 It is obvious that $(2)$ is
 dual to $(1)$. Let $r$ be the rank of $E$ and $s$ the rank of
$S$. Without loss of generality, we can assume,  at a fixed point
$p\in X$, there exists a local holomorphic frame
$\{e_1,\cdots,e_r\}$ of $E$ centered at point $p$ such that
$\{e_1,\cdots, e_s\}$ is a local holomorphic frame of $S$. Moreover,
we can assume that $h(e_\alpha,e_\beta)(p)=\delta_{\alpha\beta},
\text{for } 1\leq \alpha, \beta \leq r.$ Hence, the curvature tensor
of $S$ at point $p$ is \beq R^S_{i\bar j
\alpha\bar\beta}=-\frac{\p^2 h_{\alpha\bar\beta}}{\p z^i\p\bar
z^j}+\sum_{\gamma=1}^{s}\frac{\p h_{\alpha\bar\gamma}}{\p
z^i}\frac{\p h_{\gamma\bar\beta}}{\p\bar z^j} \eeq where $1\leq
\alpha,\beta\leq s$. The curvature tensor of $E$ at point $p$ is
\beq R^E_{i\bar j \alpha\bar\beta}=-\frac{\p^2
h_{\alpha\bar\beta}}{\p z^i\p\bar z^j}+\sum_{\gamma=1}^{r}\frac{\p
h_{\alpha\bar\gamma}}{\p z^i}\frac{\p h_{\gamma\bar\beta}}{\p\bar
z^j} \eeq where $1\leq \alpha,\beta\leq r$. By  formula
(\ref{Nakanop}), it is easy to see that \beq R^E|_{S}-R^S=\frac{\sq
}{2\pi}\sum_{i,j}\sum_{\alpha,\beta=1}^s\left(\sum_{\gamma=s+1}^r\frac{\p
h_{\alpha\bar\gamma}}{\p z^i}\frac{\p h_{\gamma\bar\beta}}{\p\bar
z^j}\right)dz^i\wedge d\bar z^j\ts e^\alpha\ts e^\beta\label{cu}\eeq
is semi-Nakano-positive. Hence $(1)$ follows.
 \eproof \elemma

\vspace{0.2cm}

\section{Curvatures of direct image sheaves and the moduli space $\sM_g$}

\subsection{Curvature of direct image sheaves}
Let $\sX$ be a K\"ahler manifold with dimension $d+n$ and $S$ a
K\"ahler manifold with  dimension $d$. Let $f:\sX\>S$  be a proper
K\"ahler fibration. Hence, for each $s\in S$,
$$X_s:=f^{-1}(\{s\})$$
is a compact K\"ahler manifold with dimension $n$. Let
$(\sE,h^{\sE})\>\sX$ be a Hermitian holomorphic vector bundle.
Consider the space of holomorphic $\sE$-valued $(n,0)$-forms on
$X_s$,
 $$E_s:=H^0(X_s,\sE_s\ts K_{X_s})\cong H^{n,0}(X_s,\sE_s)$$
where $\sE_s=\sE|_{X_s}$. It is well-known that, if the vector
bundle $\sE$ is  ``positive" in certain sense, by Grauert locally
free theorem, there is a natural holomorphic structure on
$$E=\bigcup_{s\in S}\{s\}\times E_s$$  such that the vector bundle $E$ is isomorphic to the direct
image sheaf $f_*(K_{\sX/S}\ts \sE)$.
 Using the canonical isomorphism
$K_{\sX/S}|_{X_s}\cong K_{X_s},$
 a \emph{local smooth} section $u$ of $E$ over $S$ can be identified as
 a family of
$\sE$-valued holomorphic $(n,0)$ form on $X_s$. By this
identification, there is a natural metric on $E$. For any local
smooth section $u$ of $E$, one can define a Hermitian metric on $E$
by \beq h(u,u)=c_n\int_{X_s} \left\{u,u\right\}\label{metric}\eeq
where $c_n=(\sq)^{n^2}$. Here, we only use the Hermitian metric of
$\sE_s$ on each fiber $X_s$ and we do not specify  background
K\"ahler metrics on the fibers.

In particular, we consider a canonically polarized family
$f:\sX\>S$. We  define the Hodge metric on $E=f_*(K^{\ts
m}_{\sX/S})$ where $m$ is an integer with $m\geq 2$. Let $L_t$ be
the restriction of the line bundle $\sL:=K_{\sX/S}^{\ts (m-1)}$ on
the fiber $X_t$. Let $U\subset S$ be a small open neighborhood of
$t=0$. Let $s$ be a local smooth section of $f_*(K_{\sX/S}^{\ts
m})$,  by the very definition of direct image sheaf
$f_*(K_{\sX/S}^{\ts m})$, for any $t\in U$ \beq s(t)\in
H^0(X_t,K_{X_t}^{\ts m})\eeq By this identification, for any local
smooth sections $s_\alpha,s_\beta$ of $f_*(K_{\sX/S}^{\ts m})$, we
define the \emph{Hodge metric} on $f_*(K_{\sX/S}^{\ts m})$ by \beq
h(s_\alpha,s_\beta)=c_n\int_{X_t}\langle s_\alpha,s_\beta\rangle_{}
dV_{t}\label{Hodge}\eeq  where $\langle\bullet,\bullet\rangle_{}$ is
the pointwise inner product on $\Gamma(X_t,K^{\ts m}_{X_t})$. More
precisely, if $s_\alpha=\phi_\alpha\ts e$ and $s_\beta=\phi_\beta\ts
e$, where $e$ is a local holomorphic basis of $K^{\ts m}_{X_t}$ and
$\phi_\alpha,\phi_\beta$ are local smooth functions on $X_t$, the
metric is
$$h(s_\alpha,s_\beta)=c_n\int_{X_t}\langle s_\alpha, s_\beta\rangle_{} dV_t=c_n\int_{X_t}\la \phi_\alpha,\phi_\beta\ra |e|^2_{} dV_t$$
where $|e|^2_{}$ is the canonical metric on $K^{\ts m}_{X_t}$. Since
the family is canonically polarized, the metrics on
$E=f_*(K_{\sX/S}^{\ts m})$ defined by (\ref{metric}) and
(\ref{Hodge}) are the same.

Let $(t_1,\cdots, t_d)$ be the local holomorphic coordinates
centered at a point $p\in S$ and $$\nu: T_pS\>H^1(X,TX)$$ the
Kodaira-Spencer map on the center fiber $X=\pi^{-1}(p)$ and
$\theta_i\in\H^{0,1}(X,T^{1,0}X)$ the harmonic representatives of
the images $\nu(\frac{\p}{\p t^i})$ for $i=1,\cdots, d$. Let
$\{\sigma_\alpha\}$ be a basis of $\H^{n,0}(X,L)$ where
$L=K_{\sX/S}^{\ts (m-1)}|_X$. The following theorem is well-known
(e.g. \cite[Theorem IV]{Schumacher12} with $p=n$; for similar
formulations see also \cite{Liu-Sun-Yau04, Liu-Sun-Yau05,
Liu-Sun-Yau08, Berndtsson11,
 Sun-Yau11, Sun12,Liu-Yang}.)

\btheorem\label{main1234} If $f:\sX\>S$ is effectively parameterized
and $m\geq 2$,  at point $p$, the curvature tensor of the Hodge
metric $h$ on $f_*(K_{\sX/S}^{\ts m})$  is
\begin{eqnarray} R_{i\bar j \alpha\bar\beta}&=&\nonumber(m-1)\left((\Delta+m-1)^{-1}\left(
\theta_i\lrcorner\sigma_\alpha\right),
\theta_j\lrcorner\sigma_\beta\right)\\&+&
(m-1)\left((\Delta+1)^{-1}(\langle\theta_i,\theta_j\rangle)\cdot
 \sigma_\alpha, \sigma_\beta \right).\label{77}\end{eqnarray}
\etheorem

\bremark Note that there are two Green's operators in formula
(\ref{77}) and they have different geometric meanings. More
precisely,  $(\Delta+m-1)^{-1}$ acts on sections and
$(\Delta+1)^{-1}$ acts on functions. When $m=2$, the curvature
formula (\ref{77}) is, in fact, different from Wolpert's curvature
formula (\cite{Wolpert86}) since in Wolpert's formula, two Green's
operators are the same and both of them act on functions. We will
analyze them carefully in the next subsection.\eremark

\subsection{Wolpert's curvature formula}

In this subsection, we will derive Wolpert's curvature formula of
the Weil-Petersson metric on the moduli space $\sM_k$ of Riemann
surfaces with genus $k\geq 2$ from  curvature formula (\ref{77}).

Let $(X_0,\omega_g)$ be a Riemann surface with the Poincar\'e metric
$\omega_g$. The Weil-Petersson metric on the moduli space $\sM_k$ is
defined  as

\beq \left(\frac{\p}{\p t_i}, \frac{\p}{\p
t_j}\right)_{WP}=\int_{X_0}\theta_i\cdot\bar \theta_j dV\eeq where
$\theta_i, \theta_j\in \H^{0,1}(X_0,T^{1,0}X_0)$ are the images of
$\frac{\p}{\p t_i}, \frac{\p}{\p t_j}$ under the Kodaira-Spencer
map: $T_0\sM_k\>\H^{0,1}(X_0,T^{1,0}X_0)$ respectively. It is known
that $f_*(K^{\ts 2}_{\sT_{k}/\sM_k})$ is isomorphic to the
holomorphic cotangent bundle $T^{*1,0}\sM_k$. Hence, there are two
Hermitian metrics on this bundle, one is the Weil-Petersson metric
and the other one is the Hodge metric defined in (\ref{Hodge}). In
order to discuss the relations between these two metrics, we can
consider the natural isomorphism
$$T:\Om^{1,0}(X_0,K_{X_0})\>\Om^{0,1}(X_0,T^{1,0}X_0)$$ given by
\beq  T\left(\eta d z\ts dz\right)=g^{-1}\bar \eta d\bar z\ts
\frac{\p}{\p z}.\eeq

\blemma\bd\item The operator $T$ is well-defined;

\item $\sigma\in \H^{1,0}(X_0,K_{X_0})$ if and only if $T(\sigma)\in
\H^{0,1}(X_0,T^{1,0}X_0)$. \ed \elemma

\noindent Let $\sigma_\alpha\in \H^{1,0}(X_0,K_{X_0})$. To simplify
notations, $T(\sigma_\alpha)$ is denoted by $\theta_\alpha$. The
local inner product on the space $\Om^{0,1}(X_0,T^{1,0}X_0)$ is
denoted by $\langle\bullet,\bullet\rangle$ and sometimes it is also
denoted by $\cdot$. That is, if
$\eta,\mu\in\Om^{1,0}(X_0,T^{1,0}X_0)$,
$\langle\eta,\mu\rangle=\eta\cdot\bar\mu$. It is easy to see that
\beq g\cdot\langle \theta,T(\sigma) \rangle d\bar z\ts
dz=\theta\lrcorner \sigma\label{87}\eeq

 \blemma The Hodge metric
coincides with the Weil-Petersson metric on the moduli space
$\sM_k$. More precisely, $T:\left(\H^{1,0}(X_0,K_{X_0}),
g_{\mathrm{Hodge}}\right)\> \left(\H^{0,1}(X_0,T^{1,0}X_0),
g_{\mathrm{WP}}\right)$ is a conjugate-isometry. \bproof Let
$\sigma_\alpha=f_\alpha dz\ts e,\sigma_\beta=f_\beta dz\ts e$, then
$$\langle \sigma_\alpha,\sigma_\beta\rangle_{\mathrm{Hodge}}=\sq\int_{X_0}g^{-1}f_\alpha\bar f_\beta dz\wedge d\bar z=\int_{X_0}g^{-2}f_\alpha\bar f_\beta dV$$
Similarly,
$$\langle T(\sigma_\beta),T(\sigma_\alpha)\rangle_{\mathrm{WP}}=\int_{X_0}g^{-2}f_\alpha\bar f_\beta dV=\langle\bar{T(\sigma_\alpha)},\bar{T(\sigma_\beta)}\rangle_{\mathrm{WP}}.$$
That is $$\langle
\sigma_\alpha,\sigma_\beta\rangle_{\mathrm{Hodge}}=\langle\bar{T(\sigma_\alpha)},\bar{T(\sigma_\beta)}\rangle_{\mathrm{WP}}.$$
\eproof

\elemma

Let $\Delta=\bp\bp^*+\bp^*\bp$ be the Laplacian operator on the
space $\Om^{p,q}(X_0,L_0)$ and $\Delta_0=\bp^*\bp$ the Laplacian
operator on $C^\infty(X_0)$.

\blemma We have the following relation between two different Green's
operators \beq
(\Delta+1)^{-1}(\theta_i\lrcorner\sigma_\alpha)=g(\Delta_0+1)^{-1}(\theta_i\cdot\bar\theta_\alpha)d\bar
z\ts e\label{twogreen}\eeq

\bproof It is obvious that the right hand side of (\ref{twogreen})
is a well-defined tensor. Hence, without loss of generality, we can
verify formula (\ref{twogreen}) in the normal coordinate of the
K\"ahler-Einstein metric. Let $\{z\}$ be the normal coordinate
centered at a fixed point $p$, i.e.,
$$g(p)=1, \frac{\p g}{\p z}(p)=\frac{\p g}{\p \bar z}(p)=0$$
The K\"ahler-Einstein condition is equivalent to \beq\Delta_0
g=-1\label{KE}\eeq at the fixed point $p$. Let $s=fd\bar z\ts e\in
\Om^{0,1}(X_0,K_{X_0})$,  at  $p$, we have \be\Delta
s&=&\bp\bp^*s=\bp\left(\left(\bp^*(fd\bar z)\right)\ts e+
g^{-1}f\frac{\p\log g}{\p z}\ts e\right)\\&=&(\Delta_0 f)d\bar z\ts
e+fd\bar z\ts e\\
&=&\left((\Delta_0+1) f\right)d\bar z\ts e\ee where we use the
K\"ahler-Einstein condition (\ref{KE}). Hence, at point $p$, \be &&
\Delta\left(g(\Delta_0+1)^{-1}(\theta_i\cdot\bar\theta_\alpha)d\bar
z\ts e\right)\\
&=&\left((\Delta_0+1)\left(g(\Delta_0+1)^{-1}(\theta_i\cdot\bar\theta_\alpha)\right)\right)d\bar
z\ts e\\ &=&g(\theta_i\cdot\bar \theta_\alpha)d\bar z\ts
e+(\Delta_0g)\left((\Delta_0+1)^{-1}(\theta_i\cdot\bar\theta_\alpha)d\bar
z\ts
e\right)\\
&=&g(\theta_i\cdot\bar \theta_\alpha)d\bar z\ts
e-\left((\Delta_0+1)^{-1}(\theta_i\cdot\bar\theta_\alpha)d\bar z\ts
e\right)\\
&=&\theta_i\lrcorner
\sigma_\alpha-g\left((\Delta_0+1)^{-1}(\theta_i\cdot\bar\theta_\alpha)d\bar
z\ts e\right) \ee where we use (\ref{87}) and $g(p)=1$ in the last
step. That is, at the fixed point $p$, (\ref{twogreen}) holds.
\eproof

\elemma  Now we obtain the well-known Wolpert formula:
\btheorem[\cite{Wolpert86}] The curvature tensor of the
Weil-Petersson metric
 on  the \bf{cotangent bundle} of the moduli space  is:
\beq R_{i\bar j \alpha\bar\beta}=\int(\Delta_0+1)^{-1}\left(
\theta_i\cdot\bar\theta_\alpha\right)\left( \bar
\theta_j\cdot\theta_\beta\right)dV+
\int(\Delta_0+1)^{-1}(\theta_i\cdot\bar \theta_j)
 \left(\bar\theta_\alpha\cdot\theta_\beta \right)dV\label{97}\eeq
\bproof If we set $m=2$ in formula (\ref{77}),  (\ref{97}) follows
from formulas (\ref{77}), (\ref{twogreen}) and (\ref{87}). \eproof
 \etheorem

\btheorem\label{main12} The Weil-Petersson metric is
dual-Nakano-negative and semi-Nakano-negative. \bproof By duality
(e.g. Remark \ref{dual123}), we only need to prove the curvature
tensor (\ref{97}) is Nakano-positive and semi-dual-Nakano-positive.
At first, we prove the semi-dual-Nakano positive part. That is, for
any nonzero matrix $u=(u^{i\alpha})$, it suffices to show  \beq
R_{i\bar j\alpha\bar\beta}u^{i\beta}\bar u^{j\alpha}\geq
0\label{prove}\eeq Here we use similar ideas of \cite[Section
4]{Liu-Sun-Yau08}. Let $G(z,w)$ be the kernel function of the
integral operator $(\Delta_0+1)^{-1}$. It is well-known that $G$ is
strictly positive and in a neighborhood of the diagonal, $G (z, w) +
\frac{1}{2\pi}\log |z- w|$ is continuous. So we obtain \be R_{i\bar
j\alpha\bar\beta}u^{i\beta}\bar
u^{j\alpha}&=&\int_{X_0}\int_{X_0}G(z,w)\theta_i(w)\bar\theta_\alpha(w)\bar\theta_j(z)\theta_\beta(z)u^{i\beta}\bar
u^{j\alpha}dV_wdV_z
\\&+&\int_{X_0}\int_{X_0}G(z,w)\theta_i(w)\bar\theta_j(w)\bar\theta_\alpha(z)\theta_\beta(z)u^{i\beta}\bar
u^{j\alpha}dV_wdV_z\ee If we set
$H(w,z)=\theta_i(w)\theta_\beta(z)u^{i\beta}$, \be R_{i\bar
j\alpha\bar\beta}u^{i\beta}\bar
u^{j\alpha}&=&\int_{X_0}\int_{X_0}G(z,w)H(w,z)\bar{H(z,w)}dV_wdV_z\\&&+\int_{X_0}\int_{X_0}G(z,w)H(w,z)\bar{H(w,z)}dV_wdV_z\ee
Since the Green's function is symmetric, i.e., $G(z,w)=G(w,z)$,
$$R_{i\bar
j\alpha\bar\beta}u^{i\beta}\bar
u^{j\alpha}=\frac{1}{2}\int_{X_0}\int_{X_0}G(z,w)\left(H(w,z)+H(z,w)\right)\bar{\left(H(w,z)+H(z,w)\right)}dV_wdV_z
$$
which is non-negative. Hence we get (\ref{prove}).

 For the Nakano-positivity, we can use the
same method. It is easy to see that we can get strict
Nakano-positivity since the Kodaira-Spencer map is injective.
 \eproof \etheorem

\bremark In virtue of Lemma \ref{optimal}, the moduli space $\sM_k$
has the same curvature property as the unit disk with the invariant
Bergman metric. It is optimal in the sense that the curvature can
not be Nakano-negative at any point. In fact, it follows from the
$L^2$-vanishing theorems on $\sM_k$ (e.g. \cite{Ohsawa82}). \eremark

\section{Curvature properties of the moduli space of curves}

In this section, we investigate the curvature properties of the
Weil-Petersson metric as well as its background Riemannian metric on
the moduli space  of curves, based on very general curvature
relations on K\"ahler manifolds.

\subsection{Curvatures on Riemannian manifold}

Let $(M,g)$ be a Riemannian manifold with Levi-Civita connection
$\nabla$. The  curvature tensor is defined as \beq
R(X,Y)Z=\nabla_X\nabla_YZ-\nabla_Y\nabla_XZ-\nabla_{[X,Y]}Z \eeq for
any $X,Y,Z\in \Gamma(M,TM)$. In the local coordinates $\{x^i\}$ of
$M$, we adopt the convention: \beq
R(X,Y,Z,W)=g(R(X,Y)Z,W)=R_{ijk\ell}X^iY^jZ^kW^\ell .\eeq The
curvature operator is \beq \sR: \Gamma(M,\Lambda^2TM)\to
\Gamma(M,\Lambda^2TM) \qtq{and} g(\sR(X\wedge Y), Z\wedge
W)=R(X,Y,W,Z)\eeq Note here, we change the orders of $Z,W$ in the
full curvature tensor. For Riemannian sectional curvature, we use
\beq K(X,Y)=\frac{R(X,Y,Y,X)}{|X|_g^2|Y|_g^2-\la X,Y\ra_g^2} \eeq
for any linearly independent vectors $X$ and $Y$.

Let $(M,g)$ be a Riemannian manifold.  $T_\C M:=TM\ts \C$ is the
complexification of the real vector bundle $TM$. We can extend the
metric $g$ and $\nabla$ to $T_\C M$ in the $\C$-linear way and still
denote them by $g$ and $\nabla$ respectively.

\bdefinition Let $(M,g)$ be a Riemannian manifold and $R$ be the
complexified Riemmanian curvature operator. We say $(X,g)$ has
non-positive (resp. non-negative) complex sectional curvature, if
\beq R(Z,\bar W, W,\bar Z)\leq 0\ \ \ \ \ \ (\textrm{resp.} \geq
0)\label{cscdef}\eeq for any $Z,W\in T_\C M$. \edefinition

\bdefinition A vector $v\in T_\C M$ is called isotropic if
$g(v,v)=0$. A subspace is called isotropic if every vector in it is
isotropic. $(M,g)$ is called to have non-positive (resp.
non-negative) isotropic curvature if \beq g\left(\sR(v\wedge w),
v\wedge w\right)\leq 0\ ( \text{resp}. \geq 0),\eeq for every pair
of vectors $v, w\in T_\C M$ which span an isotropic $2$-plane.
\edefinition

\subsection{Curvature relations on K\"ahler manifolds}

\noindent In  \cite{Siu80}, Siu introduced the following
terminology: \bdefinition Let $(X,g)$ be a compact K\"ahler
manifold. $(X,g)$ has  strongly negative curvature(resp. strongly
positive) if   \beq   R_{i\bar j k\bar \ell}\left(A^i\bar
B^j-C^i\bar D^j\right)\bar{\left(A^\ell\bar B^k-C^\ell\bar D^k
\right)}\leq 0 \ \ \text{(resp. $\geq 0$)} \label{siu}\eeq  for any
$A=A^i\frac{\p}{\p z^i}$, $B=B^j\frac{\p}{\p z^j}$,
$C=C^i\frac{\p}{\p z^i}$, $D=D^j\frac{\p}{\p z^j}$ and the identity
in the above inequality holds if and only if $ A^i\bar B^j-C^i\bar
D^j=0$ for any $i,j$.
 \edefinition

\btheorem\label{equivalent1} Let $(X,g)$ be a K\"ahler manifold.
Then $g$ is a metric with strongly non-negative curvature (resp.
strongly non-positive curvature) in the sense of Siu if and only if
the complex sectional curvature is non-negative (resp.
non-positive). \bproof Let $Z, W\in T_\C X$. In local holomorphic
coordinates $\{z^i\}$ of $X$, one can write $\displaystyle
Z=a^i\frac{\p}{\p z^i}+b^{\bar i}\frac{\p}{\p \bar z^i},\ \ \
W=c^j\frac{\p}{\p z^j}+d^{\bar j}\frac{\p}{\p \bar z^j}.$ We can
compute \begin{eqnarray} &&\nonumber R(Z,\bar W,W,\bar
Z)\\&&=R\left(a^i\frac{\p}{\p z^i}+b^{\bar i}\frac{\p}{\p \bar
z^i},\bar{c^j}\frac{\p}{\p \bar z^j}+\bar{d^{\bar j}}\frac{\p}{\p
z^j}, c^k\frac{\p}{\p z^k}+d^{\bar k}\frac{\p}{\p \bar z^k},
\bar{a^\ell}\frac{\p}{\p \bar z^\ell}+\bar{b^{\bar
\ell}}\frac{\p}{\p z^\ell}\right).\label{csc}
\end{eqnarray}It has sixteen terms, but it is well-known that on a K\"ahler
manifold $R_{ijk\ell}=0, \ R_{\bar i jk\ell}=R_{i\bar j
k\ell}=R_{ij\bar k \ell}=R_{ijk\bar\ell}=0,$ and their conjugates
are also zero, i.e. $R_{\bar i\bar j\bar k\bar \ell}=0, \ \ R_{ i
\bar j\bar k\bar \ell}=R_{\bar i j \bar k\bar \ell}=R_{\bar i\bar j
k \bar \ell}=R_{\bar i\bar j\bar k\ell}=0.$ Since $g$ is K\"ahler,
by Bianchi identity, we see $R_{ij\bar k\bar\ell}=-R_{j\bar k i\bar
\ell}-R_{\bar k i j\bar \ell}=-R_{j\bar k i\bar \ell}+R_{i\bar k
j\bar \ell}=0.$ Similarly, we have $R_{\bar i\bar j k\bar \ell}=0.$
Hence (\ref{csc})  contains  four nonzero terms, i.e., \be R(Z,\bar
W,W,\bar Z)&=&R_{i\bar j k\bar \ell}\cdot a^i\cdot\bar{c^j}\cdot
c^k\cdot \bar{a^\ell}+R_{i\bar j \bar k\ell}\cdot a^i\cdot
\bar{c^j}\cdot d^{\bar k}\cdot \bar{b^{\bar \ell}}\\&&+R_{\bar i j
\bar k \ell}\cdot b^{\bar i}\cdot \bar{d^{\bar j}}\cdot d^{\bar
k}\cdot\bar{b^{\bar \ell}}+R_{\bar i j k\bar \ell}\cdot b^{\bar i}
\cdot\bar{d^{\bar j}}\cdot c^k \cdot\bar{a^\ell}\\&=&R_{i\bar j
k\bar \ell}\cdot a^i\cdot\bar{c^j}\cdot c^k\cdot
\bar{a^\ell}-R_{i\bar j\ell \bar k}\cdot a^i\cdot \bar{c^j}\cdot
d^{\bar k}\cdot \bar{b^{\bar \ell}}\\&&+R_{j\bar i  \ell\bar k
}\cdot b^{\bar i}\cdot \bar{d^{\bar j}}\cdot d^{\bar
k}\cdot\bar{b^{\bar \ell}}-R_{j\bar i k\bar \ell}\cdot b^{\bar i}
\cdot\bar{d^{\bar
j}}\cdot c^k \cdot\bar{a^\ell}\\
&=&R_{i\bar j k\bar \ell}\cdot a^i\cdot\bar{c^j}\cdot c^k\cdot
\bar{a^\ell}-R_{i\bar j  k\bar \ell}\cdot a^i\cdot \bar{c^j}\cdot
d^{\bar \ell}\cdot \bar{b^{\bar k}}\\&&+R_{i\bar j k \bar \ell}\cdot
b^{\bar j}\cdot \bar{d^{\bar i}}\cdot d^{\bar \ell}\cdot\bar{b^{\bar
k}}-R_{i\bar j k\bar \ell}\cdot b^{\bar j} \cdot\bar{d^{\bar
i}}\cdot c^k \cdot\bar{a^\ell}\\
&=&R_{i\bar j k\bar \ell}\Big(a^i\cdot\bar{c^j}\cdot c^k\cdot
\bar{a^\ell}- a^i\cdot \bar{c^j}\cdot d^{\bar \ell}\cdot
\bar{b^{\bar k}}\\&&\ \ \ \ \ \ \ \ \ + b^{\bar j}\cdot \bar{d^{\bar
i}}\cdot d^{\bar \ell}\cdot\bar{b^{\bar k}}- b^{\bar j}
\cdot\bar{d^{\bar i}}\cdot
c^k \cdot\bar{a^\ell}\Big)\\
&=&R_{i\bar j k\bar \ell}\left(a^i\cdot\bar{c^j}-b^{\bar j}\cdot
\bar{d^{\bar i}}\right)\left(c^k\cdot \bar{a^\ell}-d^{\bar
\ell}\cdot\bar{b^{\bar k}}\right)\\
&=&R_{i\bar j k\bar \ell}\left(a^i\cdot\bar{c^j}-b^{\bar j}\cdot
\bar{d^{\bar i}}\right)\bar{\left({a^\ell}\cdot \bar{c^k}-b^{\bar k}
\cdot\bar{d^{\bar \ell}}\right)}.\ee Let $A^{i\bar
j}=a^i\cdot\bar{c^j}-b^{\bar j}\cdot \bar{d^{\bar i}}$. We obtain
\beq R(Z,\bar W, W,\bar Z)=R_{i\bar j k\bar \ell} A^{i\bar
j}\cdot\bar{A^{\ell\bar k}}=R_{i\bar j k\bar \ell} A^{i\bar
\ell}\cdot\bar{A^{j\bar k}}.\eeq Therefore, the curvature is
strongly non-negative  in the sense of Siu if and only if the
complex sectional curvature is non-negative. The proof for the
equivalent on non-positivity is similar.
 \eproof \etheorem

\btheorem\label{curvatureop} If $(X,g)$ is a K\"ahler manifold with
semi dual-Nakano-positive curvature (resp. semi dual-Nakano-negative
curvature ), then its background Riemmanian curvature operator is
non-negative (resp. non-positive). \bproof Let $z^i=x^i+\sq y^i$ be
the local holomorphic coordinates centered at a given point. Then
from the relation $\displaystyle\frac{\p }{\p
z^i}=\frac{1}{2}\left(\frac{\p}{\p x^i}-\sq \frac{\p}{\p
y^i}\right),\ \ \ \ \frac{\p }{\p \bar
z^i}=\frac{1}{2}\left(\frac{\p}{\p x^i}+\sq \frac{\p}{\p
y^i}\right)$ one obtains $\displaystyle\frac{\p}{\p
x^i}=\frac{\p}{\p z^i}+\frac{\p}{\p \bar z^i},\ \ \ \ \ \frac{\p}{\p
y^i}=\sq \left(\frac{\p}{\p z^i}-\frac{\p}{\p\bar z^i}\right).$ On
the background Riemannian manifold,  any vector $V$ in $\Lambda^2
T_\R X$ can be written as \beq V=a^{ij}\frac{\p}{\p x^i}\wedge
\frac{\p}{\p x^j}+b^{pq}\frac{\p}{\p x^p}\wedge \frac{\p}{\p y^q}+
c^{mn}\frac{\p}{\p y^m}\wedge \frac{\p}{\p y^n}.\eeq In the
coordinates $\{z^i,\bar z^i\}$, we have \be
V&=&a^{ij}\left(\frac{\p}{\p z^i}+\frac{\p}{\p \bar
z^i}\right)\left(\frac{\p}{\p z^j}+\frac{\p}{\p \bar z^j}\right)+\sq
b^{pq}\left(\frac{\p}{\p z^p} +\frac{\p}{\p \bar
z^p}\right)\left(\frac{\p}{\p z^q}-\frac{\p}{\p\bar
z^q}\right)\\&&-c^{mn}\left(\frac{\p}{\p z^m}-\frac{\p}{\p\bar
z^m}\right)\left(\frac{\p}{\p z^n}-\frac{\p}{\p\bar z^n}\right)\\
&=& a^{ij}\left(\frac{\p}{\p z^i}\wedge \frac{\p}{\p
z^j}+\frac{\p}{\p\bar z^i}\wedge \frac{\p}{\p z^j}+\frac{\p}{\p
z^i}\wedge \frac{\p}{\p\bar z^j}+\frac{\p}{\p\bar z^i}\wedge
\frac{\p}{\p\bar z^j}\right)\\
&&+\sq b^{pq}\left(\frac{\p}{\p z^p}\wedge \frac{\p}{\p
z^q}+\frac{\p}{\p\bar z^p}\wedge \frac{\p}{\p z^q}-\frac{\p}{\p
z^p}\wedge \frac{\p}{\p\bar z^q}-\frac{\p}{\p\bar z^p}\wedge
\frac{\p}{\p\bar z^q}\right)\\
&&-c^{mn}\left(\frac{\p}{\p z^m}\wedge \frac{\p}{\p
z^m}-\frac{\p}{\p\bar z^m}\wedge \frac{\p}{\p z^n}-\frac{\p}{\p
z^m}\wedge \frac{\p}{\p\bar z^n}+\frac{\p}{\p\bar z^m}\wedge
\frac{\p}{\p\bar z^n}\right)\\
&=&A^{ij} \frac{\p}{\p z^i}\wedge \frac{\p}{\p z^j}+B^{i\bar j}
\frac{\p}{\p z^i}\wedge \frac{\p}{\p \bar z^j}+C^{\bar i\bar
j}\frac{\p}{\p \bar z^i}\wedge \frac{\p}{\p \bar z^j}\ee where
$$A^{ij}:=a^{ij}+\sq b^{ij}-c^{ij},\ \ \ \ C^{\bar i\bar j}:=a^{ij}-\sq b^{ij}-c^{ij}$$
and
$$B^{i\bar j}:=a^{ij}-\sq b^{ij}+c^{ij}-a^{ji}-\sq b^{ji}-c^{ji}.$$
By the elementary facts that
$$R_{ij\ell k}=R\left(\frac{\p}{\p z^i}\wedge \frac{\p}{\p z^j}, \frac{\p}{\p z^k}\wedge \frac{\p}{\p z^\ell}\right)=0;$$
$$R_{ij\bar \ell\bar k}=R\left(\frac{\p}{\p z^i}\wedge \frac{\p}{\p z^j}, \frac{\p}{\p \bar z^k}\wedge \frac{\p}{\p \bar z^\ell}\right)=0;$$
$$R_{i\bar j \ell k}=R\left(\frac{\p}{\p z^i}\wedge \frac{\p}{\p \bar z^j}, \frac{\p}{\p z^k}\wedge \frac{\p}{\p z^\ell}\right)=0;$$
and also their conjugates are all zero, we obtain \be
\sR(V,V)&=&R\left(B^{i\bar j} \frac{\p}{\p z^i}\wedge \frac{\p}{\p
\bar z^j}, B^{k\bar \ell} \frac{\p}{\p z^k}\wedge \frac{\p}{\p \bar
z^\ell}\right)\\&=&R_{i\bar j \bar \ell k}B^{i\bar j}B^{k\bar
\ell}=-R_{i\bar j k\bar \ell}B^{i\bar j}B^{k\bar \ell}.\ee Let
$E^{i\bar j}:=a^{ij}+c^{ij}-a^{ji}-c^{ji}, F^{i\bar
j}:=-b^{ij}-b^{ji},$ then $B^{i\bar j}=E^{i\bar j}+\sq F^{i\bar j}.$
Note that the matrix  $(E^{i\bar j})$ is real and  skew-symmetric;
the matrix $(F^{i\bar j})$ is real and symmetric. Hence, by the
curvature property \be \sR(V,V)&=&-R_{i\bar j k\bar \ell}B^{i\bar
j}B^{k\bar
\ell}\\
&=&-R_{i\bar j k\bar \ell}E^{i\bar j}E^{k\bar \ell}-\sq R_{i\bar j
k\bar\ell}(E^{i\bar j} F^{k\bar \ell}+E^{k\bar \ell} F^{i\bar
j})+R_{i\bar j k\bar \ell} F^{i\bar j}F^{k\bar \ell}.\ee On the
other hand $R_{i\bar j k\bar \ell}$ is skew-symmetric in the pairs
$(i, j)$ and  $(k,\ell)$, we obtain
$$\sq R_{i\bar j
k\bar\ell}(E^{i\bar j} F^{k\bar \ell}+E^{k\bar \ell} F^{i\bar
j})=R_{i\bar j k\bar \ell} F^{i\bar j}F^{k\bar \ell}=0$$ since
$(E^{i\bar j})$ is real and skew-symmetric and $(F^{i\bar j})$ is
real and symmetric.
 Therefore,
$$\sR(V,V)=-R_{i\bar j k\bar \ell}E^{i\bar j}E^{k\bar \ell}=-R_{i\bar j k\bar \ell} E^{i\bar \ell} E^{k\bar j}=R_{i\bar j k\bar \ell} E^{i\bar \ell} \bar{E^{j \bar k}}$$
where in the last step we use again the fact that $E^{k\bar j}$ is
real and skew-symmetric, i.e. $E^{k\bar \ell}=-{E^{j \bar
k}}=-\bar{E^{j \bar k}}$ . Now, we see that if $(X,g)$ is semi-
dual-Nakano-positive (resp. semi dual-Nakano-negative), then the
Riemannian curvature operator is non-negative (resp. non-positive).
 \eproof \etheorem

\bremark $(\P^2, \omega_{FS})$  is dual-Nakano-positive, but the
Riemannian curvature operator of the background Riemannian metric is
only non-negative. In fact, on  any compact K\"ahler manifold, there
does not exist a Riemannian metric with quasi-positive Riemannian
curvature operator since it has nonzero second Betti number (e.g.
\cite[p.212]{Petersen06}).
 \eremark

\emph{The proof of Theorem \ref{1}}. $(1)\Longrightarrow(2)$ follows
from Theorem \ref{curvatureop}, and $(3)\Longleftrightarrow(4)$
follows from Theorem \ref{equivalent1}. $(2)\Longrightarrow (4)$:
let $Z,W\in T_\C M$. Let $Z\wedge \bar W=V+iU$, where $V$ and $U$
are real tensors. Then \be R(Z,\bar W, W,\bar Z)&=&\sR(Z\wedge \bar
W, \bar
Z\wedge W)\\
&=&\sR(V+iU,V-iU)\\
&=&\sR(V,V)+\sR(U,U) \ee where the last step follows since our
curvature operator is extended to $T_\C M$ in the $\C$-linear way
and $\sR(U,V)=\sR(V,U)$. The other relations follow from similar
computations.

\bremark \bd
\item Exactly the same relations hold for semi-positivity.

\item  There is another notion called ``weakly $\frac{1}{4}$-pinched negative Riemannian sectional
curvature". If $(X,\omega)$ is a compact K\"ahler manifold with
weakly $\frac{1}{4}$-pinched negative Riemannian sectional
curvature, then $(X,\omega)$ is semi dual-Nakano-negative. Indeed,
Yau-Zheng
 proved in \cite{Yau-Zheng91}(see also \cite{Hernandez91}) that any
 {compact}
K\"ahler manifold with  weakly $\frac{1}{4}$-pinched negative
Riemannian sectional curvature must be a ball quotient. However,
Mostow-Siu surfaces (\cite{Mostow-Siu80}) have dual-Nakano-negative
curvature tensors, but they are not covered by a $2$-ball. \ed
\eremark

\emph{The proof of Theorem \ref{3}}. It follows from Theorem \ref{1}
and Theorem \ref{2}.

\section{Totally geodesic submanifolds in Torelli locus.}

In this section, we study the existence of certain locally symmetric
submanifold in moduli space  $\sM_k$ of curves with genus $k\geq 2$
by using the curvature properties we obtained. As an application of
Theorem \ref{3} and Lemma \ref{monotone}, we derive

\bcorollary\label{submonotone} Let $S$ be any submanifold of
$(\sM_k, \omega_{\mathrm{WP}})$ with the induced metric, then $S$
has

\bd \item semi Nakano-negative curvature;

\item strictly negative holomorphic bisectional curvature.

\ed \ecorollary

We need the following rigidity result by W-K. To:

\btheorem[\cite{To89}]\label{To}  Let $(X,g)$ be a locally symmetric
Hermitian manifold of finite volume uniformized by an irreducible
bounded symmetric domain of rank $\geq 2$. Suppose $h$ is Hermitian
metric on $X$ such that $(X,h)$ carries non-positive holomorphic
bisectional curvature. Then $h=cg$ for some constant $c>0$.
\etheorem

 \emph{The proof of Theorem \ref{4}}. It is well-known that any totally geodesic submanifold of
$\sA_k$ is also locally  symmetric. Suppose $j(V)$ has rank$>1$.
Then by Siu's computation in \cite[Appendix, Theorem~4]{Siu80a}, the
holomorphic bisectional curvatures of the canonical metrics on
irreducible bounded symmetric domains of rank$>1$  are non-positive
but \emph{not strictly negative}. Let $h$ be the metric on $j(V)$
induced by the Hodge metric on $\sA_k$. Hence, $(j(V), h)$ has
non-positive holomorphic bisectional curvature by formula
(\ref{cu}). Let $g$ be the Hermitian metric on $V$ induced by the
Weil-Petersson metric on $\sM_k$. By Corollary \ref{submonotone},
$(V,g)$ has strictly negative holomorphic bisectional curvature.
Since the Torelli map $j$ is holomorphic and injective, by Theorem
\ref{To}, $h=cg$ for some positive constant $c$ which is a
contradiction. Hence, $j(V)$ is of rank $1$, i.e. a ball quotient.\\

\emph{The proof of Corollary \ref{5}}. Suppose $\Om$
 has rank $>$ $1$. Then by a result of Mok ( \cite[Theorem
4]{Mok07} or \cite{Mok87}), we see $X$ must be totally geodesic in
$(\sM_k,\omega_{\mathrm{WP}})$. That is, the second fundamental form
of the immersion must be zero. Since $(\sM_k, \omega_{\mathrm{WP}}
)$ is dual-Nakano-negative, we see $(X, h)$ is also
dual-Nakano-negative, and in particular, $(X,h)$ has strictly
negative curvature in the sense of Siu which is a
contradiction(\cite[Appendix, Theorem~4]{Siu80a}). Hence $\Om$ must
be of rank $1$, i.e. $X$ must be a ball quotient.

\end{document}